\documentclass[12pt]{amsart}
\usepackage{amssymb}

\numberwithin{equation}{section}
\theoremstyle{plain}
\newtheorem{theorem}{Theorem}[section]

\newtheorem{cor}[theorem]{Corollary}

\theoremstyle{definition}

\theoremstyle{remark}

\def\Cpx{{\mathbb C}}

\def\Rn{{{\mathbb R}^n}}

   \newcommand{\Omegak}{{\Omega^{(k)}}}
  \newcommand{\Osing}{{\Omega^{{\rm sing}}}}
  \newcommand{\rank}{{\; \rm rank \;}}
  \newcommand{\spann}{{\; \rm span \;}}

\title{An open problem in complex analytic geometry
arising in harmonic analysis}

\author[Michael Ruzhansky]{Michael Ruzhansky}
\address{
  Michael Ruzhansky:
  \endgraf
  Department of Mathematics
  \endgraf
  Imperial College London
  \endgraf
  180 Queen's Gate, London SW7 2AZ 
  \endgraf
  United Kingdom
  \endgraf
  {\it E-mail address} {\rm m.ruzhansky@imperial.ac.uk}
  }
\thanks{The paper was written while the
 author was supported by the JSPS Invitational
 Research Fellowship.}

\date{\today}
\begin{document}
\begin{abstract}  In this paper, an open problem
in the multidimensional complex analysis is presented 
that arises
in the investigation of the regularity properties of 
Fourier integral operators and in the regularity theory for
hyperbolic partial differential equations. The problem 
is discussed in a self-contained elementary way and some 
results towards its resolution are presented. A conjecture 
concerning the structure of appearing affine fibrations is
formulated.
\end{abstract}

\maketitle

\section{Introduction}

In this paper, we present an open problem that arises in
harmonic analysis related to Fourier integral operators
and hyperbolic partial differential equations.
In \cite{SSS}, Seeger, Sogge, and Stein formulated
a so-called smooth factorization 
condition that guarantees the local $L^p$ bounds
for non-degenerate Fourier integral operators. 
In \cite{Ru2} it was shown that this condition is
satisfied in a number of important cases from
the point of view of the theory of hyperbolic
partial differential equations. 
Moreover, an approach to this condition
based on the notion of so-called
affine fibrations has been developed in \cite{Ru2}.
In \cite{RuP}, regularity properties of Fourier integral operators
for certain parameter dependent affine fibrations have been 
established. The partial understanding of this
problem from the point of view of harmonic analysis
can already be efficiently applied to the regularity theory
of Fourier integral operators with real (\cite{Ru2}) and
complex (\cite{Ru3}) phase functions.

In this paper we will reformulate this condition
entirely in the language of complex
analytic geometry of several complex variables. 
Subsequently, we will present results that have
been established for this problem by the author.
This analysis together with available examples suggests
that the structure of the set of singularities of 
affine fibrations must have a rather rigid nature.
In this paper we will formulate the corresponding conjecture 
that is motivated by results and examples presented here.
Combined with other methods developed in \cite{Ru2}, the
validity of this conjecture would imply the sharp
Sobolev $L^p$ properties of the corresponding
non-degenerate Fourier integral operators
and of solutions to the Cauchy problems for 
corresponding hyperbolic
partial differential equations.

In the following section we formulate the problem,
present some results and examples that motivate
a conjecture about the set of essential singularities
of the corresponding set of affine fibrations.
 In the last section we will explain how this problem
   arises naturally in the symplectic geometry and the
   subsequent theory of Fourier integral operators.

\section{Formulation of the problem}

In this section we will formulate the problem in question in
a simple self-contained way in the language of the complex
analysis of several variables.

   Let $\Omega$ be an open connected subset of
   $\Cpx^n$, and let 
   $\Gamma:\Omega\to\Cpx^n$ be a holomorphic mapping
   (in general, the image space does not have to have the same
   dimension as $\Omega$, but we consider this case for
   simplicity). 
   Let $D\Gamma(\xi)\in\Cpx^{n\times n}$ denote the
   Jacobi matrix of first order derivatives of $\Gamma$, so that
   its entries are
   $\{D\Gamma(\xi)\}_{jk}=\partial_{\xi_k}\Gamma_j(\xi),$
   $j,k=1,\ldots,n$, where we denote
   $\Gamma=(\Gamma_1,\ldots,\Gamma_n).$ 
   
   Let us make the following assumptions about the mapping $\Gamma$:
   \begin{itemize}
   \item[(A1)] 
   Let us denote
   $$k:=\max_{\xi\in\Omega}\rank D\Gamma(\xi)$$ 
   and assume that $1\leq k\leq n-1$.  
   Let $\Omegak$ denote
   the set of $\xi$ where the rank of $D\Gamma$ is maximal, i.e.
   $$\Omegak=\{\xi\in\Omega: \rank D\Gamma(\xi)=k\}.$$
   
   \item[(A2)] Assume that for every $\xi\in\Omegak$, 
     the level set 
     $$\Gamma^{-1}(\Gamma(\xi))=\{\eta\in\Omega: \Gamma(\eta)=
     \Gamma(\xi)\}$$
     is an affine $(n-k)$--dimensional space through $\xi$.
   \end{itemize}
   
   Let us now discuss these conditions. In fact, condition (A1)
   is not restrictive since it excludes the 
   cases $k=0$ and 
   $k=n$ only, where the subsequent problem is trivial.
   The case of $k=0$ would mean that the rank of $D\Gamma$ is
   identically zero, which would mean that $\Gamma$ is
   a constant mapping. The case of $k=n$ would mean that
   the level sets of $\Gamma$ are just points. Since we
   will be interested in the geometry of the level sets
   of $\Gamma$, we exclude these two cases from the consideration.
   
   Condition (A2) is the main condition that we impose on
   $\Gamma$. In fact, for $\xi\in\Omegak$, by the implicit
   function theorem we may conclude that the level set 
   $\Gamma^{-1}(\Gamma(\xi))$ is a smooth analytic 
   $(n-k)$-dimensional submanifold of $\Omega$. Condition
   (A2) assumes that all these level sets are affine
   (i.e. they are linear spaces with the origin at $\xi$),
   thus imposing a rigid structure on the geometry of the
   problem.
   
   We can also note that the level sets 
   $\Gamma^{-1}(\Gamma(\xi))$ can not intersect in $\Omegak$
   (unless they coincide).
   Indeed, if they would intersect, the intersection point
   would belong to two level sets. Thus, the mapping 
   $\Gamma$ would take the same value on both level sets.
   Since both of them are
   linear by (A2), and if the intersection point is in $\Omegak$,
   the joint level set can be regular only if they coincide.
   In this way the set $\Omegak$ becomes a union of
   non-intersecting affine spaces which are the level
   sets of $\Gamma$. So, we can write
   $$\Omegak=\bigcup_{\xi\in\Omegak}\left(
   \Gamma^{-1}(\Gamma(\xi))\cap\Omegak\right),$$
   where any two affine spaces $\Gamma^{-1}(\Gamma(\xi))$
   from the union are either 
   disjoint in $\Omegak$ or coincide.
   Thus, we will talk about an {\em affine fibration} in
   $\Omegak$, given by the union of disjoint planes
   (which are level sets of a holomorphic mapping $\Gamma$).
   
   It follows from (A1) that the set $\Omegak$
   is open and dense in $\Omega$. 
   Let us now denote
   $$\varkappa:\xi\mapsto\ker D\Gamma(\xi),\;\;\; 
   \varkappa:\Omegak\to{\mathbb G}_{n-k}(\Cpx^n),$$
   so that by the implicit function theorem
   $\varkappa$ is a regular mapping from $\Omegak$
   to the Grassmanian ${\mathbb G}_{n-k}(\Cpx^n)$, which is
   the set of all 
   $(n-k)$--dimensional linear subspaces of $\Cpx^n$.
   We note that condition (A2) implies that for
   $\xi\in\Omegak$ we have 
   $$\Gamma^{-1}(\Gamma(\xi))=\xi+\varkappa(\xi),$$
   so that the dependence of the level sets on
   $\xi\in\Omegak$ is analytic. 
   
   The main question of this paper is when (and whether)
   this fibration extends analytically from $\Omegak$
   to the set $\Omega$. In fact, it can be shown that the
   mapping $\varkappa$ extends to a 
   meromorphic mapping in $\Omega$, so by
   the analytic graph theorem this extension is analytic if and
   only if it is continuous (see e.g. \cite{Loja}). 
   We denote this extension by $\overline{\varkappa}$.
    
   To study the extension properties of $\varkappa$ from
   $\Omegak$ to $\Omega$, let us introduce the set 
   $$\Osing=\{\xi\in\Omega:\varkappa 
   \textrm{ is not continuously extendible over } \xi \},$$
   which is the set of essential singularities of 
   $\overline\varkappa$. 
   The condition that $\varkappa$ extends regularly from
   $\Omegak$ to $\Omega$ would be equivalent to the
   condition that the set $\Osing$ is empty.

   Now we will present some results, examples, and formulate
   a conjecture concerning the structure of the set
   $\Osing$ for the mappings
   $\Gamma$ satisfying conditions (A1) and (A2).
   It turns out that it is very important to know as much
   as possible about properties of the set $\Osing$.
   The following theorems have been established in
   \cite{Ru0}, and then extended in
   \cite{Ru2}. Thus, it can be shown that the mapping
   $$\overline{\varkappa}:\Omega\to{\mathbb G}_{n-k}(\Cpx^n)$$ 
   is meromorphic in the sense of \cite{Rem}
   (see also \cite{Loja}). There, a
   mapping $\tau:X\to Y$ between complex manifolds
   $X$ and $Y$ is called meromorphic if the following
   three conditions hold:
    \begin{itemize}
   \item[(1)] For every $x\in X$ the image set $\tau(x)\subset Y$ is
          non-empty and compact in $Y$.
   \item[(2)] The graph of the mapping  $\tau$, that is the set all pairs
          $(x,y)\in X\times Y$ such that $y\in\tau(x)$, is a connected
          complex analytic subset of $X\times Y$ of dimension equal
          to the dimension of $X$.
   \item[(3)] There exist a dense subset $X^*$ of $X$, such that
          for every $x\in X^*$ the image set $\tau(x)$ consists of a
          single point.
       \end{itemize}
   This immediately implies that the set $\Osing$ of
   essential singularities of a meromorphic mapping
   $\overline\varkappa$ is analytic and we have the estimate
   $\dim\Osing\leq n-2.$ With additional analysis, we have more:
          
    \begin{theorem} \label{THM:main}
     Suppose $\Osing\not=\emptyset$. 
     Then $\Osing$ is analytic and for every
     $\xi\in\Omega$ we have the estimate
     $$\max\{k-1,n-k+1\}\leq \dim_\xi\Osing\leq n-2.$$
    Moreover, let $\xi\in\Osing$ be a regular point
    of the analytic set~$\Osing$. 
    Let $\xi=\lim_{j\to\infty}\xi_j$ be a limit of
    some sequence $\xi_j\in\Omegak$, 
    and let $\varkappa(\xi_j)\to\varkappa\in
    {\mathbb G}_{n-k}(\Cpx^n)$. 
    Then we have the inclusion $\varkappa\subset T_\xi\Osing$. 
    \end{theorem} 
   We note that a sequence $\varkappa(\xi_j)$ always has a 
   convergent subsequence since the Grassmanian is compact,
   so we can start with a sequence $\xi_j\in\Osing$ in
   Theorem \ref{THM:main} for which the corresponding sequence
   $\varkappa(\xi_j)$ converges to some 
   $\varkappa\in {\mathbb G}_{n-k}(\Cpx^n).$
  As a corollary of estimates of Theorem \ref{THM:main}, we obtain
  that affine fibrations always have regular
  extensions from $\Omegak$ to $\Omega$ in lower dimensions:
   \begin{cor} \label{cor:cor}
   We have the following statements:
   \begin{itemize}
   \item[(1)] if $n\leq 3$, then $\Osing$ is empty.
   \item[(2)] if $k\leq 2$, then $\Osing$ is empty.
   \end{itemize}
   \end{cor}
   
   We can also show that the estimates on the dimension in
   Theorem \ref{THM:main} are sharp:
   \begin{theorem}[Sharpness] \label{THM:sharp}
   Let $k$ and $d$ satisfy
   $$ 3\leq k\leq n-1 \textrm{ and } 
   \max\{k-1,n-k+1\}\leq d\leq n-2.$$
    Then there exists a holomorphic mapping
    $\Gamma:\Cpx^n\to\Cpx^n$ 
    satisfying conditions {\rm (A1), (A2)} with $k$ as above, 
    such that $\dim\Osing=d$. 
   \end{theorem}
  
   Corollary \ref{cor:cor} 
   shows that the smallest dimension for which a 
   singular
   fibration can exist is $n=4$, with $k=3$. 
   Let us show that in $\Cpx^4$ there is a fibration by lines,
   which also satisfies additional property 
   \begin{equation}\label{EQ:grad-type}
    \Gamma=\nabla\psi,
   \end{equation}
   for some holomorphic function $\psi:\Omega\to\Cpx$.
   We may call the corresponding fibrations the
   {\em fibrations of gradient type.}
   In fact, in the following section we will show how such
   fibrations naturally arise in the symplectic geometry of
   Lagrangian manifolds and related Fourier integral operators,
   in which cases condition \eqref{EQ:grad-type} is satisfied
   with a function $\psi$ obtained from function $\phi$ from the
   following section  by factoring out the homogeneous direction.
   Under \eqref{EQ:grad-type} we obviously have
   $D\Gamma=D^2\psi$, which is a symmetric matrix.
   
   Now we will give some examples of mappings in Theorem
   \ref{THM:sharp}.
   Define
   $$ \psi(\xi_1,\xi_2,\xi_3,\xi_4)=\xi_1\xi_2^2+(\xi_3-\xi_2\xi_4)^2.$$
   The Hessian of the function $\psi$ has the form
   $$D^2\psi(\xi) = \left( \begin{array}{cccc} 0 & 2\xi_2 & 0 & 0 \\
     2\xi_2 & 2\xi_1+2\xi_4^2 & -2\xi_4
    & 4\xi_2\xi_4-2\xi_3 \\
      0 & -2\xi_4 & 2
      & -2\xi_2 \\
   0 &  4\xi_2\xi_4-2\xi_3
   & -2\xi_2 & 2\xi_2^2
  \end{array} \right),$$
  with the maximal rank in property (A1) being $k=3$. Moreover, 
  $$\rank D^2\psi|_{\xi_2=\xi_3=0}=2, \;
  \rank D^2\psi|_{\xi_1=\xi_2=\xi_3=0}=1.$$
  For $\xi_2\not=0$ the
  kernel of the matrix $D^2\psi(\xi)$ is one dimensional:
  $$ \ker D^2\psi(\xi)=\spann \langle 
      \left( \frac{\xi_3}{\xi_2}-\xi_4,
      \; 0, \; \xi_2, \; 1\right)\rangle.$$
  Therefore, the mapping $\varkappa$ of the direction of
  the line $\ker D^2\psi$ corresponds to
  $\frac{\xi_3}{\xi_2}$ (after an analytic change of variables),
  and it has the essential singularities in the set
  \begin{equation}\label{EQ:exosing}
   \Osing=\{\xi_2=\xi_3=0\}.
  \end{equation}
  Let us now consider other similar families of examples.
  Functions
   $$ \psi(\xi_1,\xi_2,\xi_3,\xi_4)=
   \xi_1\xi_2^k+(\xi_3-\xi_2\xi_4)^m$$
  with $k,m\geq 2$ again 
  lead to fibrations with essential singularities
  at $\xi_2=\xi_3=0$. The fibers are given by
  $$ \ker D^2\psi(\xi)= \spann
    \langle \left( \frac{m}{k}\frac{(\xi_3-\xi_2\xi_4)^{m-1}}{\xi_2^{k-1}}
      , \; 0, \; \xi_2, \; 1\right)\rangle,$$
   and so the set of essential singularities is also given by
   \eqref{EQ:exosing}. 
  
    In $n$-dimensional space 
   $\xi=(x_1,\ldots,x_{n-3},y,z,w)$ define
    $$ \psi(\xi)=y^2\sum_{i=1}^{n-3}x_i+(z-yw)^2.$$
  The level sets of the gradient $\nabla\psi$ have dimension 
  $n-3$, and hence we have
  $k=3$. One can check that $\dim\Osing=n-2$
  (this also follows from dimension estimates of Theorem 
   \ref{THM:main} with $k=3$).

  Finally, consider
   $$\psi(x,y,z,v,w,s,t)=xy^2+sv^2+(z-yw-vt)^2.$$
  The maximal rank of the Hessian $D^2\psi$ equals $k=5$ and
  its kernel is spanned by the vectors
  $$\left(\frac{z-yw-vt}{y},0,0,y,1,0,0\right), 
  \left(0,0,0,v,0,\frac{z-yw-vt}{v},1\right).$$
  For $y=0$ and $z-vt=0$ the first vector has an essential
  singularity, while the second vector is continuous at
  $v\not=0$. At the same time, for $v=0$ and
  $z-yw=0$, the second vector has an essential
  singularity while the first one is continuous for
  $y\not=0$. Hence
  $\{y=0, z=vt\}\cup\{v=0,z=yw\}\subset\Osing$, 
  and thus $\dim\Osing=n-2$ (by the dimension estimates of 
  Theorem \ref{THM:main}). The intersection of spaces
  $\{y=0, z=vt\}$ and $\{v=0,z=yw\}$ is 4-dimensional.

 In all
 the examples we observe that the set $\Osing$ has an
 some affine structure: it is either affine itself, or
 it is a union of affine spaces. In fact, since the condition
 (A2) imposes a lot of affine structure for the whole problem,
 it may be the case that $\Osing$ inherits these properties
 as well. So, we may formulate the following conjecture,
 as least in the case when $\Gamma=\nabla\psi$ for some
 holomorphic function $\psi$.
 
 \bigskip
   \noindent
   {\bf Conjecture.}  
  {\em  Let $\Omega$ be an open connected subset of $\Cpx^n$ and let
  $\psi:\Omega\to\Cpx$ be a holomorphic function. Suppose that
  the mapping $\Gamma=\nabla\psi$ satisfies properties}
  (A1), (A2). {\em Then the set
  $\Osing$ of essential singularities of the corresponding
  fibration $\varkappa$ is a union of affine spaces.}
 
 \bigskip
 We also note that 
 these examples present singular fibrations for different
 dimensions of the fibers. On the other hand, the dimension
 of $\Osing$ in all examples equals $n-2$. Under certain
 conditions it can be proved using dimension estimates in
 Theorem \ref{THM:main}, or otherwise, that $\dim_\xi\Osing=n-2$
 at all (regular) points.
 However, our constructions are based on
 the same idea for the singularity. It would be
 interesting to investigate whether the condition 
 $\dim\Osing=n-2$ is necessary for fibrations 
 satisfying \eqref{EQ:grad-type}, or whether all affine
 spaces in the conjecture above would be $(n-2)$-dimensional.
    
  \section{Relation to the symplectic geometry and 
     Fourier integral operators} 
   
  The problem of the analysis of the set
  $\Osing$ under conditions (A1) and (A2) 
  is closely related to the 
  theory of Fourier integral operators
  and the microlocal analysis in the following way. 
  
  Let $M$ be a smooth real analytic manifold. 
  Let $T^*M$ denote the cotangent bundle
  of $M$, and let $\pi:T^*M\to M$ be the canonical projection.
  Let $\Lambda\subset T^*M$ be a conic analytic Lagrangian
  submanifold of $T^*M\backslash 0$ endowed
  with the standard symplectic form. Since we assume that all the
  spaces are analytic, we may always extend them to the complex
  domain, so it does not matter much whether we formulate
  the following relation in real or complex language.
    
  Let $\Sigma$ be the regular part of $\pi(\Lambda)$,
  and we denote the conormal bundle of $\Sigma$ by
  \begin{equation}\label{EQ:conormal} 
    N^*\Sigma=\{ (m,\zeta)\in T^*M: m\in\Sigma, \zeta(\delta m)=0, 
     \forall \delta m\in T_m\Sigma\}.
  \end{equation}
  Then the 
  canonical projection $\pi$ defines the following 
  restrictions: 
  $$ \begin{array}{ccccc}
     N^*\Sigma & \subset & \Lambda & \subset & T^*M \\
     \downarrow &     & \downarrow &       & \downarrow \pi \\
     \Sigma & \subset & \pi(\Lambda) & \subset & M
     \end{array}
  $$
  It turns out that the generating function for the
  restriction $\pi: N^*\Sigma\to \Sigma$ leads to
  a mapping satisfying
  conditions (A1) and (A2). Let us explain this in
  more detail.

  By the Poincar\'e lemma of the symplectic
  geometry (see e.g. \cite{Duis}), we can locally 
  parameterise the
  Lagrangian manifold $\Lambda$ by some
  smooth analytic function $\phi$. This means that
  at a point $(m,\zeta)\in\Lambda$, we can locally write
  $\Lambda$ as
     $$\Lambda=\{(\nabla\phi(\zeta),\zeta)\}.$$
  It turns out that the fibers of 
  the mapping $d\pi|_\Lambda$
  correspond precisely to the level sets of $\ker D^2\phi$,
  which are linear spaces in $\zeta$ in the conormal bundle
  \eqref{EQ:conormal}.
  Thus, if we define  $$\Gamma=\nabla\phi,$$ then
  conditions (A1) and (A2) are satisfied for $\Gamma$.
  We can reduce the dimension further by using the fact that
  $\Lambda$ is conic, so we can factor out the conic
  variable from $\phi$, to obtain another function
  $\psi$, which would define $\Gamma$ by \eqref{EQ:grad-type}.
  
  Consequently, the set $\Osing$ corresponds to the
  set of points where the fibration of $N^*\Sigma$ 
  (as a union of affine spaces) is not continuously extendible to 
  $\Lambda$. The details of this construction can be found
  in \cite{Ru2} or in \cite{Ru3}.
  
  It was shown in \cite{Ru2} that 
  the set $\Osing\cap\Rn$
  is precisely
  the set where the so-called smooth factorization condition
  of Seeger, Sogge, and Stein \cite{SSS} breaks down.

  The smooth factorization condition is the sufficient 
  condition for non-degenerate Fourier integral operators
  (of suitable orders) to be locally bounded in $L^p$. 
  In turn, solutions to the
  Cauchy problems for strictly hyperbolic 
  partial differential equations can be written as a 
  sum of Fourier integral operators applied to the
  Cauchy data (see e.g. \cite{Duis}), so the problem is
  also directly related to the local $L^p$-estimates for
  solutions to hyperbolic equations (see \cite{Ru2} for
  details). 
  
  Let us briefly review now how the problem
  in terms of the symplectic geometry above is related to
  Fourier integral operators.
  Let $X$ and $Y$ be smooth paracompact real (analytic)
  manifolds of dimension $n$. Let $\Lambda^\prime$
  be a conic (analytic) Lagrangian manifold in 
  $(T^*X\backslash0)\times (T^*Y\backslash 0)$, 
  equipped with the symplectic form
  $\sigma_X\oplus -\sigma_Y$.
 Let $T\in I^\mu(X,Y;\Lambda^\prime)$ be a Fourier integral operator
 of order $\mu$ with the canonical relation
 $\Lambda^\prime$, so that
    $$\Lambda=\{(x,\xi,y,\eta): (x,\xi,y,-\eta)\in\Lambda^\prime\}
      \subset T^*(X\times Y)$$
  is the wave front set of the
  integral kernel of the Fourier integral operator
  $T$. We consider non-degenerate
  Fourier integral operators $T$,
  for which the canonical relation $\Lambda^\prime$ locally is the
  graph of a symplectomorphism between $T^*X\backslash 0$ and
  $T^*Y\backslash 0$, equipped with canonical symplectic forms
  $\sigma_X$ and $\sigma_Y$.
   By the equivalence-of-phase-function theorem
   the set $\Lambda$ locally has the form
    $$ \Lambda=\{(x,\nabla_x\Phi,y,\nabla_y\Phi):
       \nabla_\theta\Phi(x,y,\theta)=0 \},$$
   for some phase function $\Phi=\Phi(x,y,\theta$, 
   positively homogeneous
   of order one in $\theta$. With this phase,
   the Fourier integral operator $T$ can be written in the
   usual microlocal form
 \begin{equation}\label{EQ:fioT}
   Tu(x)=\int_Y\int_\Rn e^{i\Phi(x,y,\theta)}
   a(x,y,\theta) u(y) d\theta dy, 
 \end{equation}
 with some amplitude $a\in S^{\mu}_{1,0}$, a symbol of order $\mu$.
 Note, that the homogeneity of the canonical relation implies
 that $\rank d\pi_{X\times Y}|_\Lambda\leq 2n-1$, where
 $\pi_{X\times Y}$ is the canonical projection from
 $T^*(X\times Y)$ to $X\times Y$.
 The result of \cite{SSS} states
 that Fourier integral operators $T\in I^\mu(X,Y;\Lambda^\prime)$
 are locally continuous in $L^p$, provided
 $1<p<\infty$ and $\mu\leq -(n-1)|1/p-1/2|$.
  
  However, if the rank of $d\pi_{X\times Y}|_\Lambda$ does not
  attain $2n-1$, the order $\mu$ is not sharp and it depends on 
  properties of the projection
  $\pi_{X\times Y}|_\Lambda$.
  For example, it turns out that 
  the regularity properties of Fourier integral
  operators in $L^p$ spaces with $p\not=2$ depend on
  the maximal rank of the projection $\pi_{X\times Y}$ restricted to
  $\Lambda$. An important ingredient is the following {\em smooth
  factorization condition} for $\pi_{X\times Y}$.
  Suppose that there exists a number $k$, $0\leq k\leq n-1$, such that
  for every
  $\lambda_0=(x_0,\xi_0,y_0,\eta_0)\in\Lambda$, there exist a conic
  neighborhood $U_{\lambda_0}\subset\Lambda$ of $\lambda_0$ and
  a smooth homogeneous of zero order map
  $\pi_{\lambda_0}:\Lambda\cap U_{\lambda_0}\to\Lambda$ with
  constant rank,  $\rank d\pi_{\lambda_0}=n+k$, for which holds
 \begin{equation}
   \pi_{X\times Y}=\pi_{X\times Y}\circ\pi_{\lambda_0}.
   \label{eq:FIO-factor}
 \end{equation}
 Under this assumption it was shown in \cite{SSS} that
 operators $T\in I^\mu_\rho(X,Y;\Lambda^\prime)$ are continuous from
   $L^p_{comp}$ to $L^p_{loc}$, provided
   that $\mu\leq -(k+(n-k)(1-\rho))|1/p-1/2|$, with
  $1<p<\infty$ and $1/2\leq\rho\leq 1$. Denoting
  $M=X\times Y$, we obtain that condition \eqref{eq:FIO-factor}
  is equivalent to the condition that the fibration
  of $N^*(\Sigma)$ by affine fibres is smoothly extendible
  to $\Lambda$ (modulo a correction for zero sections). 
  Here $\pi=\pi_{X\times Y}$, and
  $\Sigma=\pi(N^*(\Sigma))$ is the
  regular part of $\pi(\Lambda)$, which is the singular
  support of the integral kernel of operator $T$ in
  \eqref{EQ:fioT}.


\end{document}